\documentclass[12pt,reqno]{amsart}

\usepackage{amsmath,amsfonts,amsthm,amssymb}
\usepackage{braket}

\usepackage[a4paper,height=20cm,top=24mm,bottom=24mm,left=26mm,right=26mm]{geometry}

\usepackage[numbers,sort]{natbib}
\usepackage[breaklinks=true]{hyperref}

\usepackage{graphicx}
\usepackage{epstopdf}
\usepackage{booktabs}
\usepackage[table]{xcolor}
                                        
\DeclareMathAlphabet{\mathpzc}{OT1}{pzc}{m}{it}
\usepackage{titlesec}
\titleformat{\section}
{\bfseries\scshape\centering}
{\thesection.}{.5em}{}
\titleformat{\subsection}
{\rmfamily\bfseries}
{\thesubsection}{.5em}{}
\titleformat{\subsubsection}
{\rmfamily\bfseries}
{\thesubsubsection}{.5em}{}
\numberwithin{equation}{section}

\newcommand{\I}{\mathrm{i}}
\newcommand{\E}{\mathrm{e}}

\DeclareMathOperator{\Li}{Li_2}

\DeclareMathOperator{\Vol}{Vol}

\DeclareMathDelimiter{\Norm}{\mathord}{largesymbols}{"3E}{largesymbols}{"3E}

\DeclareMathOperator{\Ad}{Ad}
\begin{document}
\abovedisplayskip=2pt
\belowdisplayskip=2pt
\baselineskip 16pt
\parskip 8pt
\sloppy


\title{Braiding Operator via Quantum Cluster Algebra}


 \author[K. Hikami]{Kazuhiro \textsc{Hikami}}

 \address{Faculty of Mathematics,
   Kyushu University,
   Fukuoka 819-0395, Japan.}

 \email{KHikami@gmail.com}

  \author[R. Inoue]{Rei \textsc{Inoue}}

  \address{Department of Mathematics and Informatics,
    Faculty of Science,
    Chiba University,
    Chiba 263-8522, Japan.}

  \email{reiiy@math.s.chiba-u.ac.jp}


\date{August 30, 2014}

\begin{abstract}
  We construct a  braiding operator in terms of the quantum dilogarithm
  function based on the  quantum cluster algebra.
  We show that
  it  is a $q$-deformation of the $\mathsf{R}$-operator
  for which hyperbolic octahedron is assigned.
  Also shown is  that,
  by taking $q$ to be a root of unity,
  our braiding operator 
  reduces to the Kashaev $\mathbf{R}^{\text{K}}$-matrix up to 
  a simple gauge-transformation.
\end{abstract}



\subjclass[2000]{}


\maketitle

\section{Introduction}

It is conjectured~\cite{Kasha96b} that the hyperbolic volume of knot
complement is given by
\begin{equation}
  \lim_{N\to\infty}
  \frac{2\pi}{N} \log
  \left| \langle K \rangle_N \right|
  = \Vol(S^3 \setminus K)
  ,
\end{equation}
where $\langle K \rangle_N$ is the Kashaev invariant of knot $K$.
Kashaev constructed  the quantum
invariant $\langle K \rangle_N$ based on the finite-dimensional
representation of the quantum
dilogarithm~\cite{Kasha95}.
It was later
realized that the Kashaev invariant $\langle K \rangle_N$ coincides
with the $N$-colored Jones polynomial at
the $N$-th root of unity~\cite{MuraMura99a},
\begin{equation}
  \langle K \rangle_N
  =
  J_N(K; q=\E^{2\pi\I/N}) ,
\end{equation}
where the $N$-colored Jones polynomial is normalized to be
$J_N(\text{unknot};q)=1$.
More precisely it was shown that Kashaev's braiding matrix~
$\mathbf{R}^{\text{K}}$
as a finite-dimensional representation of the Artin braid group
with $n$~strands,
\begin{equation}
  \label{braid_group}
  {B}_n
  =
  \left\langle
    \sigma_1, \sigma_2, \dots, \sigma_{n-1}
    ~\middle|~
    \begin{aligned}
      & \sigma_i \, \sigma_{i+1} \, \sigma_i
      =
      \sigma_{i+1} \, \sigma_{i} \, \sigma_{i+1} ,
      \\
      &\sigma_i \, \sigma_j
      =
      \sigma_j \, \sigma_i,
      \quad
      \text{for $|i-j|>1$}
    \end{aligned}
  \right\rangle ,
\end{equation}
is gauge-equivalent to
the $\mathbf{R}^{\text{J}}$-matrix for
the $N$-colored Jones polynomial at the root of unity.
As the  colored Jones polynomial is well-understood in the framework of
the quantum group
$\mathcal{U}_q(s\ell_2)$ (see, \emph{e.g.}, \cite{KassRossTura97a}),
the  Kashaev invariant $\langle K \rangle_N$
is regarded as an invariant for
$\mathcal{U}_{\E^{\pi \I/N}}(s\ell_2)$.
In~\cite{Kasha94a,Kasha96a,BaseiBened05a}
(see also~\cite{CostanMuraka13a})
studied is a
relationship between
$\langle K \rangle_N$ and
the $6j$-symbol of $\mathcal{U}_{\E^{\pi \I/N}}(s\ell_2)$,
but
the mathematical background of the Kashaev $\mathbf{R}^{\text{K}}$-matrix itself
still remains unclear, at
least, to us.

Meanwhile,
studies on the geometrical  content of $\langle K \rangle_N$ have been
much developed.
It is now recognized
that an ideal hyperbolic  octahedron is assigned to each
$\mathbf{R}^{\text{K}}$-matrix~\cite{DThurs99a},
and 
proposed~\cite{ChoMurYok09a,YYokota11a}
was
a method to 
construct from a set of  such octahedra
the Neumann--Zagier potential function~\cite{NeumZagi85a}
which give the complex volume of knot complement.
This observation is based on a fact that the hyperbolic volume of
ideal tetrahedron is given in terms of the dilogarithm function
(see, \emph{e.g.},~\cite{WPThurs80Lecture}),
and that the $\mathbf{R}^{\text{K}}$-matrix asymptotically consists 
of four dilogarithm functions~\cite{DThurs99a}.

In our previous paper~\cite{HikamiRInoue13a},
we constructed the $\mathsf{R}$-operator  from the viewpoint of the
cluster algebra.
We  showed that
the $\mathsf{R}$-operator is
geometrically interpreted as a hyperbolic octahedron which is the same
assigned to the Kashaev~$\mathbf{R}^{\text{K}}$-matrix.
The cluster algebra was  originally introduced by S.~Fomin and
A.~Zelevinsky~\cite{FominZelev02a} to study the total positivity in
semi-simple Lie groups, and
it
is promising to clarify
a deep connection with geometry~\cite{FomiShapThur08a}
(see also~\cite{NagaTeraYama11a,HikamiRInoue12a}).

The purpose of this article
is to quantize the $\mathsf{R}$-operator in~\cite{HikamiRInoue13a}.
Basic tool is the
quantum cluster
algebra~\cite{BerenZelev05a,Keller11a,FockGonc09a,FockGonc09b}.
We shall construct the  $\mathcal{R}$-operator in terms of the quantum
dilogarithm function as a conjugation for the quantum $\mathsf{R}^q$-operator, and
clarify a relationship with the Kashaev $\mathbf{R}^{\text{K}}$-matrix.
We show explicitly
that the $\mathcal{R}$-operator
reduces to the $\mathbf{R}^{\text{K}}$-matrix up
to a  gauge-transformation
when
the quantized parameter~$q^2$ tends to the~$N$-th root of unity.

This paper is organized as follows.
In Section~\ref{sec:classical}, we  briefly review our previous
results~\cite{HikamiRInoue12a,HikamiRInoue13a}.
We discuss a  relationship between the cluster algebra and the
hyperbolic geometry, and 
we recall a definition of
the $\mathsf{R}$-operator~\eqref{classicalR}
which is illustrated as a hyperbolic octahedron.
In Section~\ref{sec:quantum} we study a $q$-deformation
of the
$\mathsf{R}$-operator.
We construct the braiding $\mathsf{R}^q$-operator~\eqref{q_R_mu}
as a conjugation of the $\mathcal{R}$-operator~\eqref{R_dilog},
which is written in terms of the
quantum dilogarithm function.
In a limit that~$q^2$ goes to a root of unity, the
$\mathcal{R}$-operator reduces to the Kashaev
$\mathbf{R}^{\text{K}}$-matrix.

\section{Cluster Algebra and Hyperbolic Geometry}
\label{sec:classical}
\subsection{Cluster Algebra}

We briefly collect a notion of the cluster algebra.
See~\cite{FominZelev02a} for detail.

Fix a positive integer $N$.
Let $(\boldsymbol{x}, \mathbf{B})$ be a cluster seed,  
where~$\boldsymbol{x}=(x_1, x_2, \dots, x_N)$ is  a cluster variable,
and an $N\times N$ skew-symmetric integral matrix $\mathbf{B}=(b_{ij})$
is an exchange matrix.
The exchange matrix is depicted as a quiver which has $N$ vertices,
by regarding
\begin{equation}
  \label{b_and_arrow}
  b_{ij}
  =
  \# \left\{ \text{arrows from $i$ to $j$} \right\}
  -
  \# \left\{ \text{arrows from $j$ to $i$} \right\} .
\end{equation}


What is important is an operation on cluster seeds, which is called
the mutation. 
For $k = 1,\ldots,N$, 
the mutation $\mu_k$ of  $(\boldsymbol{x}, \mathbf{B})$ is defined by
\begin{equation}
  \mu_k(\boldsymbol{x}, \mathbf{B})
  =
  (\widetilde{\boldsymbol{x}}, \widetilde{\mathbf{B}}) ,
\end{equation}
where a cluster variable~$\widetilde{\boldsymbol{x}}=(\widetilde{x}_1, \dots,
\widetilde{x}_N)$ and an exchange matrix
$\widetilde{\mathbf{B}}=(\widetilde{b}_{ij})$ are respectively given
by
\begin{gather}
  \label{mu_on_x}
  \widetilde{x}_i
  =
  \begin{cases}
    x_i, 
    &
    \text{for $i \neq k$,}
    \\[2ex]
    \displaystyle
    \frac{1}{x_k} \,
    \left(
      \prod_{j: b_{jk}>0}
      x_j^{~b_{jk}}
      +
      \prod_{j: b_{jk} < 0}
      x_j^{~-b_{jk}}
    \right) ,
    &
    \text{for $i = k$,}
  \end{cases}
  \\
  \label{mutation_B}
  \widetilde{b}_{ij}
  =
  \begin{cases}
    - b_{ij},
    & \text{for $i=k$ or $j=k$},
    \\[2ex]
    \displaystyle
    b_{ij}
    + \frac{
      \bigl| b_{ik} \bigr| \, b_{kj}
      +
      b_{ik} \, \bigl| b_{kj} \bigr|
    }{2},
    & \text{otherwise.}
  \end{cases}
\end{gather}
In this article,
for each seed $(\boldsymbol{x}, \mathbf{B})$
we define the $y$-variable 
$\boldsymbol{y}=(y_1, \dots, y_N)$ as
\begin{equation}
  \label{y_from_x}
  y_j = \prod_k x_k^{~ b_{kj}} .
\end{equation}
The mutation of the cluster seed induces 
the mutation of the $y$-variable,
\begin{equation}
  \label{mu_on_y}
  \mu_k( \boldsymbol{y}, \mathbf{B})
  =
  (
  \widetilde{\boldsymbol{y}},
  \widetilde{\mathbf{B}}
  ) ,
\end{equation}
where the exchange matrix $\widetilde{\mathbf{B}}$ is~\eqref{mutation_B},
and
$\widetilde{y}_j = \prod_k \widetilde{x}_k^{~\widetilde{b}_{kj}}$
is given by
\begin{equation}
  \widetilde{y}_i
  =
  \begin{cases}
    y_k^{~-1} ,
    & \text{for $i = k$,}
    \\[2ex]
    y_i \, \left( 1+y_k^{~-1} \right)^{-b_{ki}} ,
    & \text{for $i \neq k$, $b_{ki} \geq 0$,}
    \\[2ex]
    y_i \left( 1+ y_k \right)^{-b_{ki}} ,
    & \text{for $i \neq k$, $b_{ki} \leq 0$.}
  \end{cases}
\end{equation}


\begin{figure}[tbph]
  \centering
  \includegraphics[]{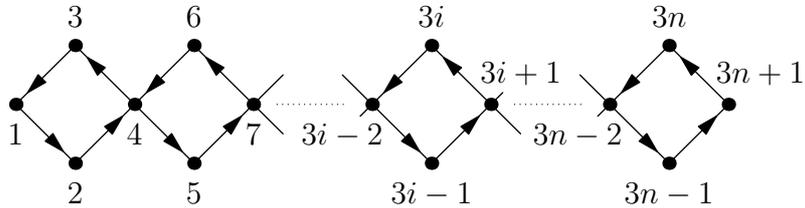}
  \caption{Quiver for the exchange matrix $\mathbf{B}$~\eqref{our_B}.}
  \label{fig:quiverR}
\end{figure}

\begin{figure}[tbph]
  \centering
  \includegraphics[]{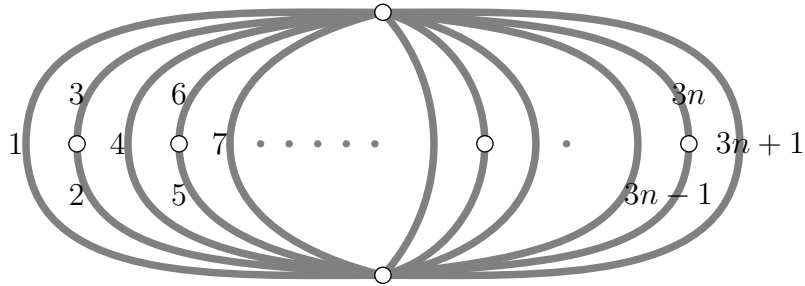}
  \caption{Triangulation of an $(n+2)$-punctured disk.
    Each edge corresponds to the vertex of the quiver in
    Fig.~\ref{fig:quiverR}.
  }
  \label{fig:triangleR}
\end{figure}

\subsection{Braiding Operator}

To study the braid group~$B_n$~\eqref{braid_group},
we set the exchange matrix   $\mathbf{B}$
to be a  $(3n+1)\times(3n+1)$ skew-symmetric matrix~\cite{HikamiRInoue13a}
\begin{align}
  \label{our_B}
  \mathbf{B} = 
  \begin{pmatrix}
    0 & 1 & -1 & 0 & 0 & \cdots  &  \cdots & \cdots & 0
    \\
    -1 & 0 & 0 & 1 & 0 & \ddots & \cdots & \cdots & \vdots
    \\
    1 & 0 & 0 & -1 & 0 & 0 & \ddots & \cdots & \vdots
    \\
    0 & -1 & 1 & 0 & 1 & -1 & 0 & \ddots & \vdots
    \\
    0 & 0 & 0 & -1 & 0 & 0 & 1 &\ddots  & 0
    \\
    \vdots & \ddots & 0 & 1 & 0 & 0 & -1 & \ddots & 0
    \\
    \vdots & \vdots & \ddots & 0 & -1 & 1 & 0  & \ddots & 1
    \\
    \vdots & \vdots & \vdots & \ddots&\ddots&\ddots & \ddots & \ddots &-1
    \\
    0 & \cdots &\cdots&\cdots& 0 &0&-1&1&0
  \end{pmatrix} .
\end{align}
See Fig.~\ref{fig:quiverR} for the associated quiver.
We remark that 
the quiver gives rise to a triangulation of a punctured disk, whose
edges correspond to the vertices in the quiver.
See Fig.~\ref{fig:triangleR}.
We define  the $\mathsf{R}$-operator 
acting on a cluster seed~$(\boldsymbol{x}, \mathbf{B})$
by
\begin{equation}
  \label{classicalR}
  \overset{i}{\mathsf{R}}
  =
  s_{3i,3i+2} \, s_{3i-1,3i+2} \, s_{3i,3i+3} \,
  \mu_{3i+1} \, \mu_{3i-1} \, \mu_{3i+3} \, \mu_{3i+1},
\end{equation}
for $i=1,\ldots,n-1$,
where
$s_{i,j}$ is  the permutation of subscripts, \emph{e.g.},
\begin{equation*}
  s_{i,j} ( \dots, x_i, \dots, x_j , \dots)
  =
  (\dots, x_j, \dots, x_i, \dots) .
\end{equation*}
As the exchange
matrix~$\mathbf{B}$ is invariant under $\overset{i}{\mathsf{R}}$,
we  write $(\overset{i}{\mathsf{R}}(\boldsymbol{x}), \mathbf{B})$
for $\overset{i}{\mathsf{R}}(\boldsymbol{x}, \mathbf{B})$
with
\begin{align}
  \label{i-R_on_x}
  \overset{i}{\mathsf{R}}(\boldsymbol{x})
  =
  \left(
    x_1, \dots, x_{3i-3},
    \mathsf{R}_x(x_{3i-2}, \dots, x_{3i+4}),
    x_{3i+5}, \dots, x_{3n+1}
  \right) ,
\end{align}
where  we have from~\eqref{mu_on_x}
\begin{gather}
    \label{R_on_x}
  \begin{aligned}[b]
    &
    \mathsf{R}_x(x_1,x_2,\dots, x_7)
    =
    \Bigl(
    x_1,
    x_5,
    \frac{x_1 \, x_3 \, x_5 + x_3 \, x_4 \, x_5
      + x_1 \, x_2 \,  x_6}{
      x_2 \, x_4} ,
    \\
    & \qquad
    \frac{
      x_1 \, x_3 \, x_4 \, x_5+ x_3 \, x_4^{~2} \, x_5 +
      x_1 \ x_3 \, x_5 \, x_7 + x_3 \, x_4 \, x_5 \, x_7
      + x_1 \, x_2 \, x_6 \, x_7}{
      x_2 \, x_4 \, x_6},
    \\
    & \qquad \qquad
    \frac{
      x_3 \, x_4 \, x_5 + x_3 \, x_5 \, x_7 + x_2 \, x_6 \, x_7}{
      x_4 \, x_6} ,
    x_3  ,
    x_7
    \Bigr) .
  \end{aligned}
\end{gather}
By definition~\eqref{y_from_x},
the action on the $y$-variable~$\boldsymbol{y}$ is induced as
\begin{align}
  \label{i-R_on_y}
  \overset{i}{\mathsf{R}}(\boldsymbol{y})
  =
  \left(
    y_1, \dots, y_{3i-3},
    \mathsf{R}_y(y_{3i-2}, \dots, y_{3i+4}),
    y_{3i+5}, \dots, y_{3n+1}
  \right) ,
\end{align}
%
%
where~\eqref{mu_on_y} gives
\begin{gather}
  \label{R_on_y}
  \begin{aligned}[b]
    &
    \mathsf{R}_y(y_1,y_2, \dots, y_7)
    =
    \Bigl(
    y_1 \, \left(
      1 + y_2 + y_2 y_4
    \right)  ,
    \frac{ y_2 \, y_4 \, y_5 \,  y_6}{
      1+ y_2 + y_6 + y_2 \, y_6 + y_2\, y_4\, y_6} ,
    \\
    & \qquad
    \frac{1+y_2+y_6+y_2 \, y_6 + y_2 \, y_4 \, y_6}{
      y_2 \, y_4},
    \frac{y_4}{
      (1+y_2+y_2 \, y_4) \, ( 1+y_6+y_4 \, y_6)
    } ,
    \\
    & \qquad
    \frac{1+y_2+y_6+y_2 \, y_6 + y_2 \, y_4 \, y_6}{
      y_4 \, y_6} ,
    \frac{ y_2 \, y_3 \, y_4 \, y_6}{
      1+ y_2 + y_6 + y_2 \, y_6 + y_2\, y_4 \, y_6} ,
    \left(
      1+y_6+y_4 \, y_6
    \right) \, y_7
    \Bigr) .
  \end{aligned}
\end{gather}
We use $\overset{i}{\mathsf{R}}$ in both~\eqref{i-R_on_x}
and~\eqref{i-R_on_y} without confusion.

\begin{figure}[tbph]
  \centering
  \includegraphics[]{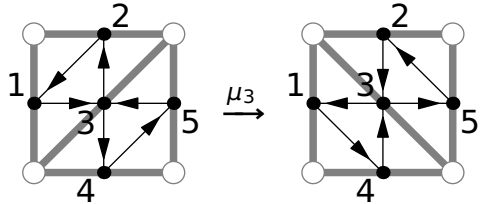}
  \caption{The mutation on a quiver is interpreted as a  flip in  triangulated
    punctured surface.
  }
  \label{fig:flip}
\end{figure}

\begin{figure}[tbph]
  \centering
  \begin{tabular}{ccc}
    \raisebox{-14mm}{\includegraphics[scale=0.8]{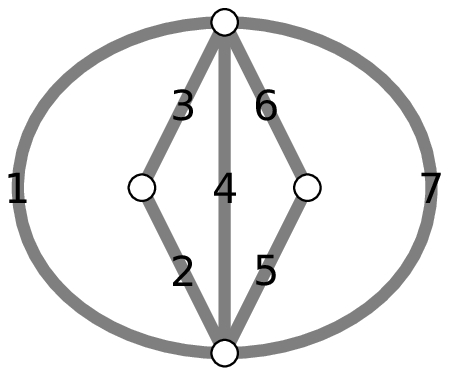}}
    $\xrightarrow[]{\mu_4}$
    &
    \raisebox{-14mm}{\includegraphics[scale=0.8]{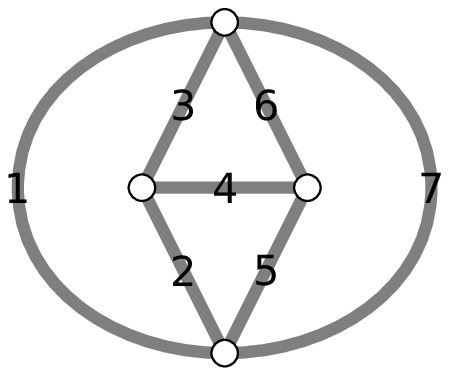}}
    $\xrightarrow[]{\mu_6}$
    &
    \raisebox{-14mm}{\includegraphics[scale=0.8]{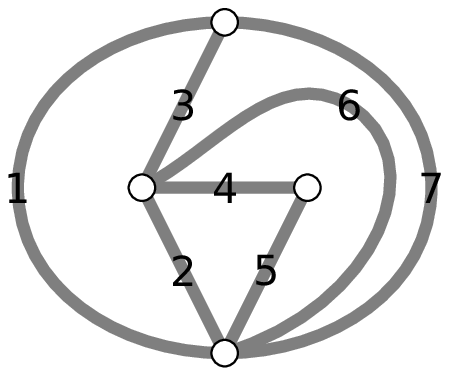}}
    \\ 
    \\
    & &
    $\phantom{\mu_2}\left\downarrow\rule{0cm}{3mm}\right. {\scriptstyle{\mu_2}}$
    \\
    \\
    \raisebox{-14mm}{\includegraphics[scale=0.8]{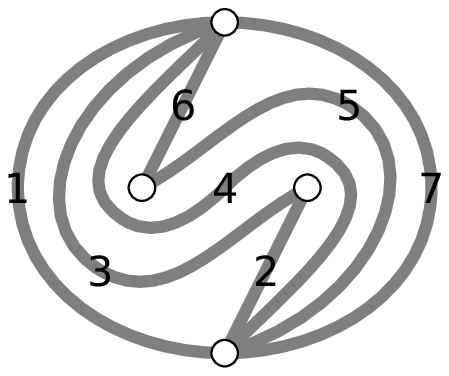}}
    $\xleftarrow[]{s}$
    &
    \raisebox{-14mm}{\includegraphics[scale=0.8]{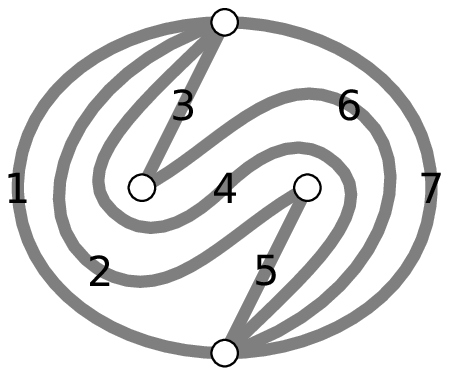}}
    $\xleftarrow[]{\mu_4}$
    &
    \raisebox{-14mm}{\includegraphics[scale=0.8]{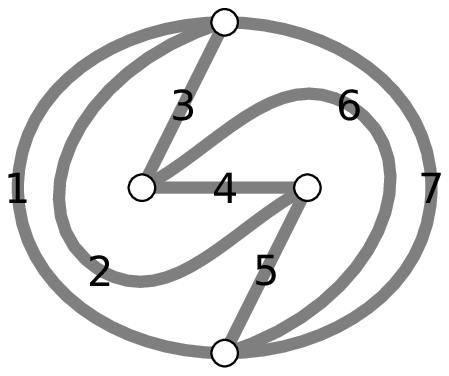}}
  \end{tabular}
  \caption{The $\mathsf{R}$-operator as a half Dehn twist.
    Illustrated is
    the action of
    $\overset{1}{\mathsf{R}}=s \, \mu_4 \mu_2 \mu_6 \mu_4$
    with $s=s_{3,5}s_{2,5}s_{3,6}$
    on the
    punctured disk.
    One sees that the $\mathsf{R}$-operator is regarded as the half
    Dehn twist exchanging two punctures counter-clockwise.    
  }
  \label{fig:R_flip}
\end{figure}

In~\cite{HikamiRInoue13a},
it was shown
that the $\mathsf{R}$-operator represents the braid group~$B_n$, and
that we have,
\begin{equation}
  \label{braid_R}
  \begin{aligned}[t]
    \overset{i}{\mathsf{R}} \, \overset{i+1}{\mathsf{R}} \, \overset{i}{\mathsf{R}}
    & =
    \overset{i+1}{\mathsf{R}} \, \overset{i}{\mathsf{R}} \,
    \overset{i+1}{\mathsf{R}}  ,
    \\
    \overset{i}{\mathsf{R}} \, \overset{j}{\mathsf{R}}
    & =
    \overset{j}{\mathsf{R}} \, \overset{i}{\mathsf{R}},
    \qquad
    \text{for $|i-j|>1$.}
  \end{aligned}
\end{equation}
This can be  proved  by direct computations using~\eqref{R_on_x}
and~\eqref{R_on_y}.
We note that 
the birational Yang--Baxter map in~\cite{Dynnik02a}
is intrinsically same with~\eqref{R_on_x}.
The braid relation~\eqref{braid_R}
could as well  be   checked  from a  dual picture as follows.
We recall that
the mutation is regarded as a ``flip'' of triangulation of
a punctured disk~\cite{FomiShapThur08a}.
Here a  flip is meant to
remove  a common edge of two adjacent triangles and to 
reproduce another  different  diagonal edge of quadrilateral
(see, for example,
Fig.~\ref{fig:flip}).
This interpretation explains
the action of the
$\overset{1}{\mathsf{R}}$-operator on a punctured disk
as
illustrated in Fig.~\ref{fig:R_flip}.
We find that the $\mathsf{R}$-operator on the punctured disk
is nothing but
a   half Dehn twist exchanging two punctures
counter-clockwise.
This  clarification of the braid group is well-known
(see, \emph{e.g.},~\cite{Birman74}),
and 
the braid relation~\eqref{braid_R} follows immediately.

\begin{figure}[tbph]
  \centering
  \begin{minipage}{0.4\linewidth}
    \centering
    \includegraphics[]{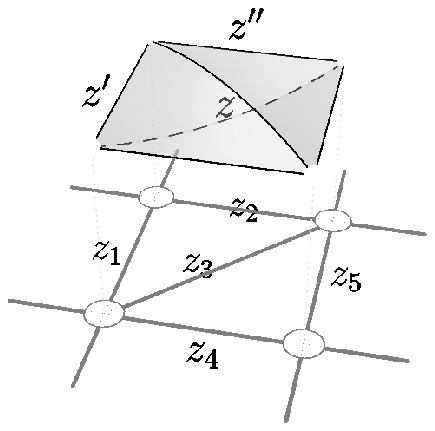}    
  \end{minipage}
  \begin{minipage}{0.4\linewidth}
    \centering
    \includegraphics[]{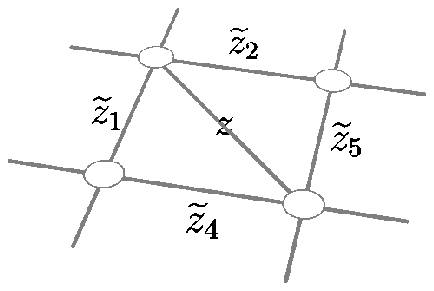}    
  \end{minipage}
  \caption{Interpretation of the flip in Fig.~\ref{fig:flip}.
    A gluing of an ideal tetrahedra to a triangulated surface
    (left)
    results in a flipped triangulated surface (right).
    Here labeled is a dihedral angle of  edge.
    Consistency condition gives~$z \, z_3=1$,
    and we have~$\widetilde{z}_1 = z_1 \, z^{\prime}$ and so on.
    This transformation is identified with the mutation of
    $y$-variable in Fig.~\ref{fig:flip}~\cite{HikamiRInoue12a}.
  }
  \label{fig:attach}
\end{figure}

\begin{figure}[tbph]
  \centering
  \includegraphics[]{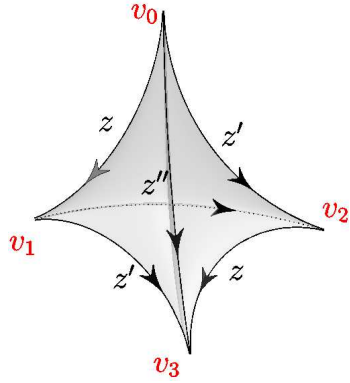}
  \caption{Hyperbolic ideal tetrahedron with a shape parameter~$z$.
    Labeled is a dihedral angle of  each edge.
    Here vertices~$v_a$ are on the
    boundary~$\partial \mathbb{H}^3$,
    and
    we mean
    $z^\prime=1-z^{-1}$ and
    $z^{\prime \prime} = (1-z)^{-1}$.
    Note that  opposite edges have the same dihedral angle.
  }
  \label{fig:tetrahedron}
\end{figure}

A geometrical interpretation of the $\mathsf{R}$-operator is  given
from the three-dimensional picture of the flip~\cite{HikamiRInoue12a}.
The mutation in Fig.~\ref{fig:flip}
acts on the $y$-variable as
\begin{align*}
  \begin{cases}
    \widetilde{y}_1 = y_1 (1+y_3) ,
    \\
    \widetilde{y}_2 = y_2 ( 1+ y_3^{-1})^{-1},
    \\
    \widetilde{y}_3 = y_3^{-1},
    \\
    \widetilde{y}_4 = y_4(1+y_3^{-1})^{-1},
    \\
    \widetilde{y}_5 = y_5 (1+y_3) .
  \end{cases}
\end{align*}
The flip in  Fig.~\ref{fig:flip} is interpreted as a
gluing  of   a
hyperbolic ideal tetrahedron  to a triangulation of punctured disk
as  in Fig.~\ref{fig:attach}.
Here
the ideal hyperbolic
tetrahedron has  a shape parameter $z$, and  
a dihedral angle of each edge is parameterized by~$z$,
$z^\prime=1-z^{-1}$, and~$z^{\prime\prime}=(1-z)^{-1}$ as in
Fig.~\ref{fig:tetrahedron}
(see, \emph{e.g.},~\cite{WPThurs80Lecture}).
Then  each dihedral angle on triangulated surface
after the gluing  is read as
\begin{equation*}
  \begin{cases}
    \widetilde{z}_1 = z_1 \, z^\prime, 
    \\
    \widetilde{z}_2 = z_2 \, z^{\prime \prime} ,
    \\
    \widetilde{z}_3 = z ,
    \\
    \widetilde{z}_4 = z_4 \, z^{\prime \prime} ,
    \\
    \widetilde{z}_5 = z_5 \, z^\prime ,
  \end{cases}
\end{equation*}
with a consistency condition $z_3 \, z=1$.
These two sets of equations indicate 
a correspondence between  the $y$-variables and the dihedral angles of
triangulated surface,
$z_k = - y_k$,
and
we conclude  that
the mutation is  regarded as a gluing  of an ideal tetrahedron with
shape parameter~$ z = -y_3^{-1}$ to
punctured surface.

\begin{figure}[tbph]
  \begin{minipage}{0.33\linewidth}
    \centering
    \includegraphics[]{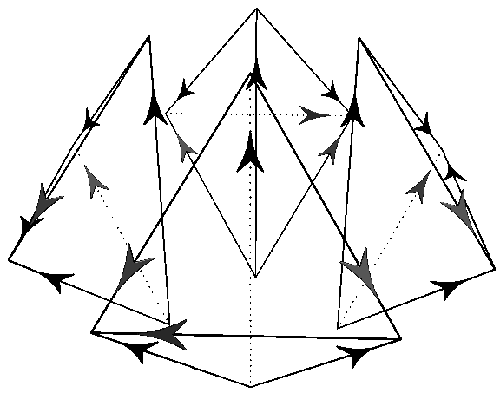}
  \end{minipage}
  \begin{minipage}{0.33\linewidth}
    \centering
    \includegraphics[]{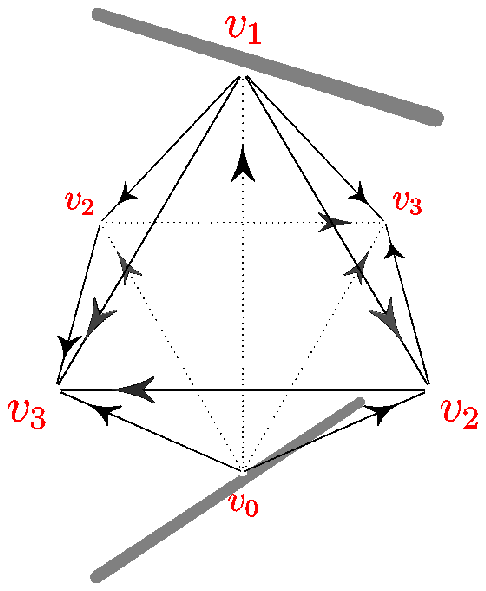}
  \end{minipage}
  \begin{minipage}{0.3\linewidth}
    \centering
    \includegraphics[]{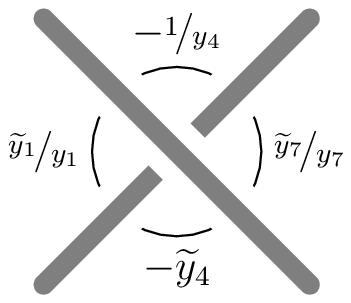}
  \end{minipage}
  \caption{Octahedron assigned to crossing (center).
    Four oriented tetrahedra (left) are glued together.
    A dihedral angle of each edge  is given in terms of the $y$-variable,
    and we  give dihedral angles around central axis
    (right).
    Here we assume
    $\widetilde{\boldsymbol{y}}=\overset{1}{\mathsf{R}}(\boldsymbol{y})$.
  }
  \label{fig:octahedron}
\end{figure}

As a consequence,
the cluster $\mathsf{R}$-operator~\eqref{classicalR}, which consists
of four mutations,
can be regarded as an ideal octahedron in Fig.~\ref{fig:octahedron}.
See that every dihedral angle is written in terms of the
$y$-variable.
Accordingly, the hyperbolic volume of
the octahedron for
$\widetilde{\boldsymbol{y}}=\overset{i}{\mathsf{R}}(\boldsymbol{y})$
is given by
\begin{equation}
  D(-1/y_{3i+1})
  +
  D(\widetilde{y}_{3i-2}/ y_{3i-2})
  +
  D(- \widetilde{y}_{3i+1})
  +
  D( \widetilde{y}_{3i+4}/ y_{3i+4}) ,
\end{equation}
where $D(z)$ is the Bloch--Wigner function
(see, \emph{e.g.},~\cite{Zagier07a}),
\begin{equation*}
  \label{BW}
  D(z)
  =
  \Im \Li(z) + \arg(1-z) \, \log | z | .
\end{equation*}
It should be noted that
this type of the hyperbolic octahedron is used not only in studies of the
Kashaev $\mathbf{R}^{\text{K}}$-matrix~\cite{DThurs99a}
but in \texttt{SnapPea} algorithm~\cite{Weeks05a}.
Note also that
the cluster variable $\boldsymbol{x}$ is identified with  Zickert's edge parameter~\cite{Zicke09}.
See~\cite{HikamiRInoue13a} for detail.

\section{Quantization}
\label{sec:quantum}

\subsection{Quantum Cluster Algebra}

We recall a quantization of the cluster algebra based
on~\cite{FockGonc09a,FockGonc09b}.
Fix a  parameter~$q$.
The $y$-variable is quantized to be a $q$-commuting generator
$\boldsymbol{Y} = (Y_1,\ldots,Y_N)$
satisfying
\begin{equation}
  \label{YYq}
  Y_k \,  Y_j = q^{2 b_{jk}} Y_j \, Y_k ,
\end{equation}
where $\mathbf{B}=(b_{ij})$ is the skew-symmetric exchange matrix used 
in the classical cluster algebra.
The $q$-commuting relation~\eqref{YYq} is realized by
\begin{equation}
  Y_k = \E^{2 \pi b \widehat{y}_k} ,
\end{equation}
where $q=\E^{\pi \I b^2}$ and 
\begin{align}\label{quantum-y}
  [\widehat{y}_k ~,~ \widehat{y}_j] = \frac{\I}{2 \pi} b_{jk}.
\end{align}
The $q$-deformation of the mutation~\eqref{mu_on_y} on the
$y$-variable
is defined by
\begin{equation}
  \label{q-mu}
  \mu^q_k(\boldsymbol{Y}, \mathbf{B})
  =
  (\widetilde{\boldsymbol{Y}},
  \widetilde{\mathbf{B}}
  ) ,
\end{equation}
where 
the exchange matrix $\widetilde{\mathbf{B}}$ is~\eqref{mutation_B},
and
\begin{equation}
  \label{quantum_Y_mutation}
  \widetilde{Y}_i
  =
  \begin{cases}
    Y_k^{-1} ,
    & \text{for $i = k$,}
    \\[2ex]
    \displaystyle
    Y_i \,
    \prod_{m=1}^{b_{ki}}
    \left( 1+q^{2m-1} \, Y_k^{-1} \right)^{-1} ,
    & \text{for $i \neq k$, $b_{ki} \geq 0$,}
    \\[2ex]
    \displaystyle
    Y_i \,
    \prod_{m=1}^{-b_{ki}}
    \left( 1+ q^{2m-1} \,  Y_k \right) ,
    & \text{for $i \neq k$, $b_{ki} \leq 0$.}
  \end{cases}
\end{equation}
One sees that this reduces to the classical mutation~\eqref{mu_on_y}
in $q\to 1$.

It is known that
the quantum mutation $\mu^q_k$~\eqref{q-mu} is decomposed into
\begin{equation}
  \mu^q_k = \mu^\sharp_k \circ \mu^\prime_k .
\end{equation}
Here $\mu^\prime_k$ is given by
\begin{equation}
  Y_i \mapsto
  \begin{cases}
    Y_k^{-1} ,
    & \text{for $i = k$,}
    \\
    q^{b_{ik} b_{ki}} \, Y_i \, Y_k^{b_{ki}} ,
    & \text{for $i \neq k$, $b_{ki} \geq 0$,}
    \\
    Y_i  ,
    & \text{for $i \neq k$, $b_{ki} \leq 0$,}
  \end{cases}
\end{equation}
and $\mu^\sharp_k$ is a conjugation by the quantum dilogarithm function
\begin{equation}
  \mu^\sharp_k
  =
  \Ad(\Phi(\widehat{y}_k)):
  Y_i \mapsto
  \Phi(\widehat{y}_k) \, Y_i \Phi(\widehat{y}_k)^{-1} .
\end{equation}
See Appendix~\ref{sec:qdilog}
for definition and properties of the quantum dilogarithm
function~$\Phi(y)$.
Note that the quantum mutation~$\mu_k^q$ is not a conjugation in general.

\subsection{Braiding Operator}

We now consider the quantum cluster algebra 
for the exchange matrix $\mathbf{B}$~\eqref{our_B}, 
whose quiver and dual picture as a triangulation
are respectively given in Fig.~\ref{fig:quiverR} and
Fig.~\ref{fig:triangleR}.

As a natural  quantization of the $\mathsf{R}$-operator~\eqref{classicalR},
we define a quantum braiding operator acting on $(\boldsymbol{Y}, \mathbf{B})$ 
by
\begin{align}
  \label{q_R_mu}
  \overset{i\phantom{q}}{\mathsf{R}^q}
  =
  s_{3i,3i+2} \, s_{3i-1,3i+2} \, s_{3i,3i+3} \,
  \mu^q_{3i+1} \, \mu^q_{3i-1} \, \mu^q_{3i+3} \, \mu^q_{3i+1} .
\end{align}
Along with the classical case 
the exchange matrix $\mathbf{B}$~\eqref{our_B} is invariant under
the operator $\overset{i\phantom{q}}{\mathsf{R}^q}$,
and~\eqref{quantum_Y_mutation} gives an action on $\boldsymbol{Y}$ as
\begin{equation}
  \label{qR_action}
  \overset{i\phantom{q}}{\mathsf{R}^q}(\boldsymbol{Y})
  =
  \left(
    Y_1, \dots, Y_{3i-3},
    \mathsf{R}_y^q(Y_{3i-2}, \dots, Y_{3i+4}),
    Y_{3i+5}, \dots, Y_{3n+1}
  \right) ,
\end{equation}
where
\begin{equation}
  \label{R-Y}
  \mathsf{R}_y^q(Y_1, \dots, Y_7)
  =
  \begin{pmatrix}
    Y_1 (1+q Y_2^{\prime}) \\
    Y_5 (1+q Y_4^{-1})^{-1}   (1+q {Y_2^{\prime}}^{-1})^{-1}
    (1+q {Y_6^{\prime}}^{-1})^{-1}
    (1+q {Y_4^{\prime\prime}}^{-1})^{-1}  \\
    {Y_2^{\prime}}^{-1} (1+q Y_4^{\prime\prime}) \\
    {Y_4^{\prime\prime}}^{-1} \\
    {Y_6^{\prime}}^{-1} (1+q Y_4^{\prime \prime}) 
    \\
    Y_3 (1+q Y_4^{-1})^{-1}
    (1+q {Y_2^{\prime}}^{-1})^{-1}   (1+q {Y_6^{\prime}}^{-1})^{-1}
    (1+q {Y_4^{\prime\prime}}^{-1})^{-1} \\
    Y_7 (1+q Y_6^{\prime}) 
  \end{pmatrix}^\top .
\end{equation}  
Here we have used
\begin{align}
  \label{eq:Y246}
  &Y_2^{\prime} = Y_2 (1 + q Y_4),
  &
  &
  Y_6^{\prime} = Y_6 (1 + q Y_4),
  &&
  Y_4^{\prime \prime}= Y_4^{-1} (1+q Y_2^{\prime})(1+q Y_6^{\prime}).
\end{align}
Clearly~\eqref{R-Y}
reduces to~\eqref{R_on_y} when $q\to 1$.

In our noncommutative algebra~\eqref{YYq} with~\eqref{our_B},
there exist central elements,
$Y_{3i-1} Y_{3i}$ ($i=1,\ldots,n$) and
$ Y_{1} Y_4 \cdots Y_{3n+1}$.
For simplicity we consider a subspace defined by
\begin{align}\label{eq:center}
  \widehat{y}_{3i-1} + \widehat{y}_{3i} = c , 
  \qquad   
  \text{for } i=1,\ldots,n,
\end{align}
where $c\in \mathbb{R}$.
In this setting,
we find that, in contrast to that~$\mu_k^q$~\eqref{q-mu}
is not an adjoint operator,
the  $\mathsf{R}^q$-operator~\eqref{q_R_mu} 
is written as a conjugation 
\begin{equation}
  \label{R_adjoint}
  \overset{i\phantom{q}}{\mathsf{R}^q}(\boldsymbol{Y})
  =
  \Ad(\overset{i}{\mathcal{R}})(\boldsymbol{Y})
  =
  \overset{i}{\mathcal{R}\vphantom{{}^{-1}}} \, \boldsymbol{Y} \,
  \overset{i\phantom{-1}}{\mathcal{R}^{-1}} ,
\end{equation} 
where
\begin{align}
  \label{R_dilog}
  \overset{i}{\mathcal{R}} 
  &=
  \Phi(\widehat{y}_{3i+1}) \, \Phi(\widehat{y}_{3i-1}) \,
  \Phi(\widehat{y}_{3i+3}) \, \Phi(\widehat{y}_{3i+1})^{-1} \,
  \theta(c+\widehat{y}_{3i+1}) .
\end{align}
See~\eqref{eq:theta-Phi} for the definition of~$\theta(z)$.
We can check~\eqref{R_adjoint}
by a direct computation using~\eqref{R_difference}.
See Appendix~\ref{sec:proof_adR}.
Furthermore we find  that
the $\mathcal{R}$-operator~\eqref{R_dilog} fulfills
the braid relation
\begin{equation}
  \label{eq:braid-1}
  \begin{aligned}[t]
    &\overset{i}{\mathcal{R}} \, \overset{i+1}{\mathcal{R}}  \,
    \overset{i}{\mathcal{R}}
    =  
    \overset{i+1}{\mathcal{R}} \,\overset{i}{\mathcal{R}} \,
    \overset{i+1}{\mathcal{R}} ,
    \\
    &\overset{i}{\mathcal{R}} \, \overset{j}{\mathcal{R}}
    = \overset{j}{\mathcal{R}} \, \overset{i}{\mathcal{R}} ,
    \qquad 
    \text{for $|i-j| > 1$.}
  \end{aligned}
\end{equation}
See Appendix~\ref{sec:proof_YBE} for the proof.
Note that the essentially same solution of the braid relation was
studied in~\cite{LFadd99b,RKasha01b}~\footnote{We thank Rinat Kashaev
  for kindly informing us.}.

As we have seen in the previous section that the classical
$\mathsf{R}$-operator~\eqref{classicalR}
on the $y$-variable  is interpreted as an   ideal
hyperbolic octahedron,
the $\mathcal{R}$-operator~\eqref{R_dilog} introduced as an adjoint
operator for the quantum $\mathsf{R}^q$-operator~\eqref{R_adjoint}
should be regarded as  a quantum content of the octahedron.
It is convincing
since the function $\Phi(z)$ reduces,
in a classical limit $b\to 0$,
to the dilogarithm function~\eqref{asymptotics_Phi},
which is related
to the hyperbolic volume
of ideal tetrahedron~\eqref{BW}.

\subsection{\mathversion{bold}Braiding Matrix at Generic $q$}

We shall give an infinite-dimensional representation of the quantum
$\mathcal{R}$-operator~\eqref{R_dilog}.
For this purpose,
we set $c=c^\prime+c^{\prime\prime}$ and
\begin{equation}
  \label{eq:y-xp}
  \begin{aligned}
    \widehat{y}_{3i-2}
    &= \widehat{x}_{i-1}-\widehat{x}_i,
    \\
    \widehat{y}_{3i-1}
    &= \widehat{p}_i + c^\prime ,
    \\
    \widehat{y}_{3i}
    &= -\widehat{p}_i + c^{\prime \prime},
  \end{aligned}
\end{equation}
where we mean $\widehat{x}_0 = \widehat{x}_{n+1} = 0$, and
$\widehat{x}_i$ and $\widehat{p}_i$ are generators of the Heisenberg
algebra,
\begin{align}
&[\widehat{x}_i ~,~ \widehat{x}_j]
= [\widehat{p}_i ~,~ \widehat{p}_j] = 0,
&&
[\widehat{x}_i ~,~ \widehat{p}_j] = \frac{\I}{2 \pi}  \delta_{ij}.
\end{align}

We define bases in coordinate and momentum spaces,
$\left| x \right\rangle$ and $\left| p \right\rangle$,
by
\begin{align*}
  &\widehat{x}_i \left| x \right\rangle = x_i \left|x\right\rangle ,
  &
  & \widehat{p}_i \left|p\right\rangle
  = p_i \left|p\right\rangle .
\end{align*}
These are orthonormal bases satisfying
\begin{gather*}
  \begin{aligned}
    &
    \left\langle x \middle|  x^\prime \right\rangle
    =
    \prod_{i=1}^n \delta(x_i - x^\prime_i),
    & \qquad
    &
    \left\langle p\middle| p^\prime\right\rangle
    = \prod_{i=1}^n \delta(p_i - p^\prime_i),
  \end{aligned}
  \\
  \left\langle x \middle| p \right\rangle
  =
  \overline{
    \left\langle p \middle| x \right\rangle
  }
  = \E^{2 \pi \I  \sum_i x_i p_i},
  \\
  \int\limits_{\mathbb{R}^n}
  \left|x\right\rangle \left\langle x \right| dx 
  =
  \int\limits_{\mathbb{R}^n}
  \left|p \right\rangle \left\langle p \right|  dp = 1,
\end{gather*}

A matrix element
$\left\langle x \right| \overset{1}{\mathcal{R}} \left| x^\prime\right\rangle$
of the $\mathcal{R}$-operator~\eqref{R_dilog}
is computed as follows.
Using~\eqref{eq:y-xp} 
we get
\begin{align*}
  & 
  \left\langle x_1, x_2 \right|
  \overset{1}{\mathcal{R}} 
  \left| x_1^\prime, x_2^\prime \right\rangle
  \\
  & 
  =
  \int\limits_{\mathbb{R}^n}
  dp
  \bra{x} \Phi(\widehat{x}_1-\widehat{x}_2) \,
  \Phi(\widehat{p}_1+c^\prime) \ket{p}  \,
  \bra{p} \Phi(-\widehat{p}_2 + c^{\prime\prime}) \,
  \Phi(\widehat{x}_1-\widehat{x}_2)^{-1} \,
  \theta(c+\widehat{x}_1-\widehat{x}_2) \ket{x^\prime}
  \\
  &
  = \Phi(x_1-x_2)
  \, \Phi(x_1^\prime-x_2^\prime)^{-1} \,
  \theta(c+x_1^\prime-x_2^\prime)
  \\
  &\qquad
  \times
  \int\limits_{\mathbb{R}} dp_1 \,
  \Phi(p_1+c^\prime) \, \E^{2 \pi \I p_1(x_1-x_1^\prime)}
  \cdot 
  \int\limits_{\mathbb{R}} dp_2  \,
  \Phi(-p_2+c^{\prime \prime}) \, 
  \E^{2 \pi \I       p_2(x_2-x_2^\prime)} .
\end{align*}
Applying~\eqref{eq:int-2} and~\eqref{eq:theta-Phi},
we obtain
\begin{multline}
  \label{R_infinite}
  \left\langle x_1, x_2 \right|
  \overset{1}{\mathcal{R}} 
  \left| x_1^\prime, x_2^\prime \right\rangle
  =
  \\
  \frac{
    \Phi(x_1-x_2) \, \Phi(x_2^\prime-x_1^\prime)
  }{
    \Phi(x_1-x_1^\prime+c_b) \,  \Phi(x_2^\prime-x_2+c_b)}
  \,
  \E^{ 2 \pi \I 
    \left(
      c_b (x_1^\prime - x_2^\prime -x_1 + x_2)
      + c^\prime(x_2^\prime -x_1)
      + c^{\prime\prime}(-x_1^\prime+x_2)
      + \frac{1}{12}(1-4 c_b^2) -
      \frac{1}{2}c^2
    \right)
  }.
\end{multline}


\subsection{\mathversion{bold}Braiding Matrix at Root of Unity 
$q^{2N} = 1$
}



In our preceding construction, we have used the quantum dilogarithm
function $\Phi(z)$~\eqref{F-q-dilog} introduced by Faddeev~\cite{LFadd95a}.
It is well known that,
due to that  $\E^{2 \pi b \widehat{y}_j}$ commute with
$\E^{2 \pi b^{-1}  \widehat{y}_k}$ for arbitrary  $j$ and $k$, 
we can replace $\Phi(z)$ by
$(-q \,  \E^{2\pi b z }; q^2)_\infty^{~-1}$,
\emph{i.e.},
we can drop a
$\overline{q}$-dependence in~\eqref{Phi_q-product}.
For simplicity, 
we
set~$c=0$ further, and 
we pay attention to
a finite-dimensional representation of the $\mathcal{R}$-operator
\begin{equation}
  \label{R_from_Y}
  \overset{1}{\mathcal{R}}
  =
  \frac{1}{
    (- q Y_4; q^2)_\infty} \cdot
  \frac{1}{
    (- q Y_2; q^2)_\infty} \cdot
  \frac{1}{
    (- q Y_6; q^2)_\infty} \cdot
  \frac{1}{
    (- q Y_4^{-1}; q^2)_\infty}  ,
\end{equation}
where we have used the $q$-Pochhammer symbol~\eqref{q-Pochhammer}.

We  rely  on a method of~\cite{BazhaReshe95}
to  construct explicitly a finite-dimensional representation of the
$\mathcal{R}$-operator~\eqref{R_from_Y}.
We set
\begin{equation}
  q= - \E^{- \frac{\varepsilon}{2 N^2}} \, \zeta^{\frac{1}{2}} ,
\end{equation}
and study a limit
$q^2\to\zeta$ by  $\varepsilon\to 0$.
Here $\zeta$ is the $N$-th root of unity, 
\begin{align}
  \label{def_zeta}
  &\zeta=\E^{ - 2\pi \I/N} ,
  &
  & \zeta^{\frac{1}{2}} = \E^{-\pi \I/N}  .
\end{align}
In a limit $\varepsilon\to 0$,
an asymptotics of the $q$-infinite product is given by~\cite{BazhaReshe95}
\begin{equation}
  (x;q^2)_\infty
  =
  \E^{- \frac{\Li(x^N)}{\varepsilon}} \,
  \sqrt{1-x^N} \,
  \prod_{k=1}^{N-1} (1 - \zeta^k x )^{-\frac{k}{N}}
  + O(\varepsilon) ,
\end{equation}
where  $\left| x \right| < 1$, and we have used the
Euler--Maclaurin formula.
We then obtain
\begin{equation}
  \label{limit_qY}
  (-q Y ; q^2)_\infty
  =
  \E^{- \frac{\Li(-Y^N)}{\varepsilon}} \, d(\zeta^{\frac{1}{2}}Y) +
  O(\varepsilon) ,
\end{equation}
where  $d(x)$ is defined by
\begin{equation}
  \label{d_func}
  d(x) = \left( 1 -x^N \right)^{\frac{N-1}{2N}} \,
  \prod_{k=1}^{N-1}
  \left(
    1 - \zeta^k x
  \right)^{- \frac{k}{N}} .
\end{equation}

We recall  that for 
$\widehat{u} \,  \widehat{v}
= \zeta \, \widehat{v} \,  \widehat{u}$
we have~\cite{BazhaReshe95}
\begin{equation}
  \label{Li_circ}
  \begin{aligned}
    \E^{- \Li(\widehat{u}^N)} * \widehat{v} 
    &=
    \widehat{v} \, \left( 1 - \widehat{u}^N  \right)^{-1/N} ,
    \\
    \E^{-\Li(\widehat{v}^N)} * \widehat{u}
    & =
    \widehat{u} \, \left( 1 - \widehat{v}^N \right)^{1/N} ,
  \end{aligned}
\end{equation}
where we mean
\begin{equation}
  \E^{\widehat{a}} * \widehat{b}
  = \lim_{\varepsilon \to 0}
  \E^{ \widehat{a}/\varepsilon} \,
  \widehat{b} \,
  \E^{- \widehat{a}/\varepsilon} .
\end{equation}
Substituting~\eqref{limit_qY} for~\eqref{R_from_Y},
we have an asymptotic behavior in $\varepsilon\to 0$
\begin{multline*}
  \overset{1}{\mathcal{R}}
  \simeq
  \E^{\frac{\Li(- Y_4^{N})}{\varepsilon}}
  \E^{\frac{\Li(- Y_2^{N})}{\varepsilon}}
  \E^{\frac{\Li(- Y_6^{N})}{\varepsilon}}
  \E^{\frac{\Li(-  Y_4^{-N})}{\varepsilon}}
  \\
  \times
  \left(
    \E^{-\Li(-Y_4^{-N})} *
    \E^{- \Li(-Y_6^N)} * \E^{-\Li(-Y_2^N)} *
    \frac{1}{d(\zeta^{\frac{1}{2}} Y_4)}
  \right) 
  \\
  \times
  \left(
    \E^{-\Li(-Y_4^{-N})} * \frac{1}{d ( \zeta^{\frac{1}{2}} Y_2 )}
  \right)
  \cdot
  \left(
    \E^{-\Li(-Y_4^{-N})} * \frac{1}{d ( \zeta^{\frac{1}{2}} Y_6 )}
  \right) \cdot
  \frac{1}{d(\zeta^{\frac{1}{2}}/Y_4)} .
\end{multline*}
With a help of~\eqref{Li_circ},
we find
\begin{equation}
  \label{asymptotics_R}
  \overset{1}{\mathcal{R}}
  \simeq
  \E^{\frac{\Li(- Y_4^{N})}{\varepsilon}}
  \E^{\frac{\Li(- Y_2^{N})}{\varepsilon}}
  \E^{\frac{\Li(- Y_6^{N})}{\varepsilon}}
  \E^{\frac{\Li(-  Y_4^{-N})}{\varepsilon}} \, \mathbf{R} .
\end{equation}
Here
we obtain  the dilogarithm factors in the right-hand side
as a  dominating term in
a  limit $\varepsilon\to 0$,
and $\mathbf{R}$ denotes the second
dominating  term given by
\begin{multline}
  \label{eq:Rmatrix}
  \mathbf{R}=
  \frac{1}{
    d(\zeta^{\frac{1}{2}} Y_4 ( 1+Y_2^N (1+Y_4^{-N}))^{\frac{1}{N}}
    (1+Y_6^N(1+Y_4^{-N}))^{\frac{1}{N}})
  } \\
  \times
  \frac{1}{
    d( \zeta^{\frac{1}{2}} Y_2 ( 1 + Y_4^{-N})^{\frac{1}{N}})
  }
  \cdot
  \frac{1}{
    d( \zeta^{\frac{1}{2}} Y_6 ( 1 + Y_4^{-N})^{\frac{1}{N}})
  } \cdot
  \frac{1}{
    d( \zeta^{\frac{1}{2}} /Y_4)
  }  .
\end{multline}
We note that  we have $q^2\to \zeta$ as $\varepsilon \to 0$, and
that
the quantum $Y$-operator in the above $\mathbf{R}$-matrix
fulfills
\begin{equation}
  \label{YY_zeta}
  Y_k \, Y_j = \zeta^{b_{jk}} Y_j \, Y_k ,
\end{equation}
with the root of unity $\zeta$~\eqref{def_zeta}.

We study a finite-dimensional
matrix representation of the second dominating term
$\mathbf{R}$
in~\eqref{asymptotics_R}.
We use
\begin{gather}
  \omega=\zeta^{-1} =\E^{2\pi\I/N},
\end{gather}
and 
we define
$w(x,y|n)$ for $n\in \mathbb{Z}_{\geq 0}$
and $x, y \in \mathbb{C}$ satisfying $x^N+y^N=1$
by
\begin{gather}
  w(x,y|0) =
  y^{\frac{1-N}{2}} \,
  \prod_{j=1}^{N-1} \left( 1 - \omega^{-j} x \right)^{\frac{j}{N}} ,
  \\
  \frac{w(x,y|n)}{w(x,y|0)}
  =
  \prod_{j=1}^n \frac{y}{1 - \omega^j x} .
  \label{def_w_n}
\end{gather}
Following a convention~\cite{BazhBaxt92b}
we  also  use  a multi-valued function of~$x$ by
\begin{equation}
  w(x,n) = w(x,\Delta(x)|n) ,
\end{equation}
where   $\Delta(x)$ is defined by
\begin{equation}
  \Delta(x) = (1-x^N)^{1/N} .  
\end{equation}
Note that
\begin{equation}
  \label{w_modulo}
  w(x,y|n+N)
  =
  w(x,y|n) ,
\end{equation}
and that
the function $w(x,n)$  is related to  the function $d(x)$ defined in~\eqref{d_func}
\begin{equation}
  \label{w_and_d}
  w(x,0)
  = \frac{1}{d(x)} .
\end{equation}
The function $w(x,y|n)$ is often used in studies of  integrable
models in
statistical mechanics, and
known are  the following identities,
\begin{gather}
  \label{w_omega_n}
  w(x,y|n) = w(\omega^n x, y| 0) ,
  \\
  \prod_{k=0}^{N-1} w(x,y|k)=1 .
\end{gather}
Furthermore we have
\begin{equation}
  \label{Fourier_w}
  \sum_{k=0}^{N-1} w(x,y|k) \, \omega^{n k} 
  =
  N \,
  \frac{
    (x/y)^{\frac{N-1}{2}}}{\lambda(y,x)} \,  
  \frac{1}{
    w\left( y,x| n -1\right)
  }
\end{equation}
See Appendix~\ref{sec:Fourier_w} for a definition of $\lambda(y,x)$
(see also~\cite{BazhBaxt92b,BazhaReshe95,Baxt91b,SergMangStro96a}).

We introduce an $N^2 \times N^2$ matrix representation of~\eqref{YY_zeta},
\begin{equation}
  \label{Y_from_XZ}
  \begin{aligned}
    Y_2 & = \omega^{\frac{1}{2}} \,
    \kappa_2 \,  \mathbf{X} \otimes \mathbf{1} ,
    \\
    Y_4 & = \omega^{\frac{1}{2}} \,
    \kappa_4 \,  \mathbf{Z} \otimes
    \mathbf{Z}^{-1} ,
    \\
    Y_6 & = \omega^{\frac{1}{2}} \,
    \kappa_6 \,
    \mathbf{1} \otimes
    \mathbf{X}^{-1} .
  \end{aligned}
\end{equation}
Here
$\kappa_a\in\mathbb{C}$, $|\kappa_a|<1$,
and
$N\times N$ matrices $\mathbf{Z}$ and $\mathbf{X}$, 
satisfying
$\mathbf{Z} \,\mathbf{ X} = \omega^{-1} \, \mathbf{X} \,  \mathbf{Z}$,
are defined by
\begin{align}
  &\left( \mathbf{Z} \right)_{j,k} = \omega^j \delta_{j,k} ,
  &
  &    \left( \mathbf{X} \right)_{j, k} = \delta_{j, k-1}  ,
\end{align}
where the  Kronecker delta has a period~$N$.
By substituting~\eqref{Y_from_XZ} for~\eqref{eq:Rmatrix},
we have
\begin{multline}
  \label{R_and_d}
  \mathbf{R}_{i j, k \ell}
  =
  \left[
    \frac{1}{
      d(\kappa_4  \, \Delta(\kappa_2^\prime)
      \, \Delta(\kappa_6^\prime) \,
      \mathbf{Z} \otimes \mathbf{Z}^{-1})}
  \right]_{i j, i j}
  \cdot
  \left[  
    \frac{1}{
      d(\kappa_2^\prime   \,  \mathbf{X} \otimes  \mathbf{1})}
  \right]_{i j , k j}
  \\
  \times
  \left[
    \frac{1}{
      d(\kappa_6^\prime \,  \mathbf{1} \otimes \mathbf{X}^{-1})}
  \right]_{k j, k \ell}
  \cdot
  \left[
    \frac{1}{
      d(\frac{1}{\omega \kappa_4}  \mathbf{Z}^{-1} \otimes  \mathbf{Z})}
  \right]_{k \ell, k \ell} ,
\end{multline}
where we have used
$\kappa_2^\prime = \kappa_2 \, \Delta(\kappa_4^{-1})$
and
$\kappa_6^\prime = \kappa_6 \, \Delta(\kappa_4^{-1})$.
By use of~\eqref{w_and_d} and~\eqref{Fourier_w}, we get
\begin{multline*}
  \mathbf{R}_{i j, k \ell}
  =
  \frac{
    \left( \kappa_2^\prime /
      \Delta( \kappa_2^\prime)
    \right)^{\frac{N-1}{2}}
  }{
    \lambda
    \left(
      \Delta(\kappa_2^\prime), 
      \kappa_2^\prime
    \right)
  }
  \,
  \frac{
    \left( \kappa_6^\prime /
      \Delta( \kappa_6^\prime )
    \right)^{\frac{N-1}{2}}
  }{
    \lambda
    \left(
      \Delta(\kappa_6^\prime), 
      \kappa_6^\prime
    \right)
  }
  \\
  \times
  \frac{
    w \left(
      \kappa_4 \Delta(\kappa_2^\prime) \Delta(\kappa_6^\prime), i-j
    \right) \,
    w \left(
      \omega^{-1}  \kappa_4^{-1}, \ell - k
    \right)
  }{
    w \left(
      \Delta(\kappa_2^\prime), i-k-1
    \right)
    \,
    w \left(
      \Delta(\kappa_6^\prime), \ell-j-1
    \right)
  } .
\end{multline*}

We set $\kappa_4=1-\delta^N$,
$0<\delta \ll 1$,
and $\kappa_2, \kappa_6>0$.
In a limit $\delta \to 0$,
we find with a help of~\eqref{def_w_n}
that a dominating term behaves as
\begin{equation}
  \label{limit_R_2}
  \mathbf{R}_{i j, k \ell}
  \propto
  \frac{(\omega)_{[i-k-1]} \,
    (\omega)_{[\ell - j -1]}
  }{
    (\omega)_{[i-j]} \,
    (\omega)_{[\ell - k -1]}
  } 
  \, \omega^{[\ell-k]}
  \cdot
  \delta^{
    N-1
      +[i-j] + [\ell - k -1]
      -[i-k-1] - [\ell -j -1]
  } .
\end{equation}
Here 
$[n] = n \mod N$ satisfying
$0 \leq [n] <N$,
and
we mean
$(\omega)_n=(\omega;\omega)_n$.
The origin of  $\omega^{[\ell - k]}$ is subtle.
It is due to that the
function~$w(x,n)=w(x,\Delta(x)|n)$ is
a multi-valued function,
$\frac{w(x,\omega y|n)}{w(x,y|n)}\propto\omega^n $,
and that
$w(\omega^{-1} \kappa_4^{-1}, \ell - k)$
which originates from~$d(\omega^{-1}\kappa_4^{-1} \mathbf{Z}^{-1}\otimes\mathbf{Z})$
in~\eqref{R_and_d}
crosses a branch-cut in getting~$\Delta(\omega^{-1}\kappa_4^{-1})$.
In~\eqref{limit_R_2}, we see that
\begin{multline*}
  \delta^{
    N-1
      +[i-j] + [\ell - k -1]
      -[i-k-1] - [\ell -j -1]
  } 
  \stackrel{\delta\to 0}{\to}
  \\
  \theta^{i,j}_{k,\ell}
  =
  \begin{cases}
    1, &
    \text{when $[i-j]+[j-\ell]+[\ell-k-1]+[k-i]=N-1$,}
    \\
    0, &
    \text{otherwise.}
  \end{cases}
\end{multline*}

As a result,
we get 
\begin{equation}
  \mathbf{R}_{i j, k \ell}
  \stackrel{\delta \to 0}{\to}
  \rho \,
  \frac{
    \omega^{-k+\ell} \,
    \theta^{i,j}_{k,\ell}
  }{
    (\omega)_{[i-j]} \,
    \overline{(\omega)}_{[k-i]} \,
    (\omega)_{[\ell-k-1]} \,
    \overline{(\omega)}_{[j-\ell]}
  } ,
\end{equation}
where $\rho$ is an irrelevant complex number.
Here we mean
$\overline{(\omega)}_n =
(\overline{\omega} ; \overline{\omega})_n$,
and
we have used an identity
\begin{equation*}
  (\omega)_{[n]} \, 
  \overline{(\omega)}_{[-n-1]} = N .
\end{equation*}
By construction, the $\mathbf{R}$-matrix fulfills the braid relation
\begin{equation}
  (\mathbf{R} \otimes 1 ) \,
  (1 \otimes \mathbf{R}) \,
  (\mathbf{R} \otimes 1 ) 
  =
  (1 \otimes \mathbf{R}) \,
    (\mathbf{R} \otimes 1 ) \,
  (1 \otimes \mathbf{R}) .
\end{equation}
One notices that this is gauge-equivalent with
the Kashaev $\mathbf{R}^{\text{K}}$-matrix~\cite{Kasha95,MuraMura99a}
\begin{align}
  (\mathbf{R}^{\text{K}})^{i,j}_{k,\ell}
  & =
  \frac{ N \, \omega^{-1+i-k} \, \theta^{i,j}_{k,\ell}}{
    (\omega)_{[i-j]} \,
    \overline{(\omega)}_{[j-\ell]}  \,
    (\omega)_{[\ell-k-1]} \,
    \overline{(\omega)}_{[k-i]}
  } .
\end{align}

To conclude, the Kashaev $\mathbf{R}^{\text{K}}$-matrix corresponds to
a  finite-dimensional representation of the
$\mathcal{R}$-operator~\eqref{R_dilog} which is
constructed based on the
quantum cluster algebra.
As we have seen that the classical
$\mathsf{R}$-operator~\eqref{classicalR} is regarded as  the
hyperbolic octahedron in Fig.~\ref{fig:octahedron} and that
a conjugation of the
$\mathcal{R}$-operator~\eqref{R_dilog}
is the quantum $\mathsf{R}^q$-operator which
reduces to the $\mathsf{R}$-operator in a limit $q\to 1$,
it is natural that both  $\mathbf{R}$- and $\mathbf{R}^{\text{K}}$-matrices 
are realized as the
octahedron in a limit $N \to \infty$.
Correspondingly
a matrix element~\eqref{R_infinite}  is an infinite-dimensional
analogue of the Kashaev $\mathbf{R}^{\text{K}}$-matrix.


\section*{Acknowledgements}
One of the authors (KH)  thanks Anatol~N. Kirillov and H.~Murakami for
communications.
RI thanks 
Y.~Terashima and M.~Yamazaki for discussions.
The work of KH is supported in part by JSPS KAKENHI Grant
Number~23340115, 24654041.
The work of RI is partially supported by JSPS KAKENHI Grant Number
22740111.

\appendix
\section{Quantum Dilogarithm}
\label{sec:qdilog}
We use the Faddeev quantum dilogarithm~$\Phi(z)$~\cite{LFadd95a} defined by
\begin{equation}\label{F-q-dilog}
  \Phi(z) =
  \exp \left(
    - \frac{1}{4} \int\limits_{\mathbb{R}+\I 0}
    \frac{
      \E^{-2 \I z w}
    }{
      \sinh(b \, w) \,
      \sinh(b^{-1} \, w)
    } \,
    \frac{d w}{w}
  \right) .
\end{equation}
Here we assume $b\in\mathbb{C}$ with $\Im b^2 >0$, and
we use
\begin{align}
  &q= \E^{\pi \I b^2 } ,
  &
  & \overline{q}=\E^{- \pi \I b^{-2}} ,
  &
  &
  c_b = \frac{\I}{2} (b+b^{-1}) .
\end{align}
It is well known that we have
\begin{equation}
  \label{Phi_q-product}
  \Phi(z)
  =
  \frac{
    (- \overline{q} \, \E^{2 \pi b^{-1}z} ; \overline{q}^2)_\infty
  }{
    ( - q \, \E^{2 \pi b z}; q^2)_\infty
  } ,
\end{equation}
where we have used the $q$-Pochhammer symbol
\begin{equation}
  \label{q-Pochhammer}
  (x;q)_n = \prod_{k=1}^n  \left( 1 - x \, q^{k-1} \right) .
\end{equation}
It is easy to see that
\begin{equation}
  \label{R_difference}
  \Phi(z \pm  \I  \,  b) 
  = (1+\E^{2 \pi b z} q^{\pm 1})^{\pm 1}  \,
  \Phi(z) .
\end{equation}
The classical dilogarithm function
is given in a limit $b\to 0$
\begin{equation}
  \label{asymptotics_Phi}
  \Phi\left( \frac{z}{2\pi b} \right)
  \sim
  \exp \left( \frac {\I}{2 \pi b^2} \, \Li(-\E^z) \right) .
\end{equation}

The most important property of~$\Phi(z)$ is
the pentagon identity~\cite{FaddKash94}
(also~\cite{Zagier07a}),
\begin{equation}
  \label{pentagon}
  \Phi(\widehat{x}) \, \Phi(\widehat{p}) 
  =
  \Phi(\widehat{p}) \, \Phi(\widehat{x} + \widehat{p}) \,
  \Phi(\widehat{x}) ,
\end{equation}
where
\begin{equation}
  \label{Heisenberg}
  \left[ \widehat{x} ~,~ \widehat{p} \right] = \frac{\I}{2 \pi} .
\end{equation}
See~\cite{KashaNakan11a} for a recent development on identities of the quantum
dilogarithm functions.
Notice that the function $\Phi(z)$ is used to construct the quantum
invariant in~\cite{KHikami06a}
(see also~\cite{DimGukLenZag09a,AnderKasha11a}).

Also
known is the Fourier transformation formula for~$\Phi(z)$ 
\begin{equation}
  \label{eq:int-2}
  \begin{aligned}[b]
    \int\limits_{\mathbb{R}}
    \Phi(z) \, \E^{2 \pi \I w z}
    \,dz
    &= \Phi(-w-c_b) \, \E^{ \I \pi w^2 - \I \pi(1-4c_b^2)/12}
    \\
    &= \frac{1}{\Phi(w+c_b)} \,
    \E^{ -2 \pi i w c_b + i
      \pi(1-4c_b^2)/12}.
  \end{aligned}
\end{equation}
See~\cite{SLWoron00a,FaddKashVolk00a,PonsoTesch00a} for detail.

We define  $\theta(z)$  by
\begin{equation}
  \label{eq:theta-Phi}
  \begin{aligned}[t]
    \theta(z)
    &= \Phi(z) \, \Phi(-z)  
    \\
    &=
    \E^{
      - \pi \I  z^2 +  \pi \I (1+2c_b^2)/6
      } .
  \end{aligned}
\end{equation}
We see  that
we have for~$\widehat{x}$ and~$\widehat{p}$
satisfying~\eqref{Heisenberg} 
\begin{align}
   \label{theta2}
   &\theta(\widehat{p})  \, \E^{2 \pi b \widehat{x}} 
   =
   \E^{2 \pi b   (\widehat{x}-\widehat{p})} \, \theta(\widehat{p}) ,
   \\
  \label{theta3}
  &\theta(\widehat{x}) \, \E^{2 \pi b \widehat{p}}
  =
  \E^{2\pi b (\widehat{x} + \widehat{p})} \, \theta(\widehat{x}) .
\end{align}

\section{Proof of~\eqref{R_adjoint}}
\label{sec:proof_adR}

We show that 
$\overset{1}{\mathcal{R}\vphantom{{}^{-1}}} \, Y_i \,
  \overset{1\phantom{-1}}{\mathcal{R}^{-1}}$ 
results in~\eqref{R-Y}.
We only give cases for $i=2, 3$ explicitly,
and
the others are obtained in a similar manner.
For a sake of  simplicity,
we write $\Phi_{\pm i} = \Phi(\pm \widehat{y}_i)$,
$\Phi_{i \pm j} = \Phi(\widehat{y}_i \pm \widehat{y}_j)$,
$\theta_{c+i} = \theta(c+ \widehat{y}_i)$,
and so on. 
For $i=2$, 
we compute as follows:
\begingroup
\allowdisplaybreaks
\begin{align*}
  &
  \overset{1}{\mathcal{R}\vphantom{{}^{-1}}}  \, Y_2 \,
  \overset{1\phantom{-1}}{\mathcal{R}^{-1}}
  \\
  &= 
  \Phi_4 \Phi_2 \Phi_6 \Phi_4^{-1} 
  \underline{\theta_{c+4} Y_2 \theta_{c+4}^{-1}} 
  \Phi_4 \Phi_6^{-1} \Phi_2^{-1} \Phi_4^{-1}
  \\
  &=
  q \,\Phi_4 \Phi_2 \Phi_6 
  \underline{\Phi_4^{-1} Y_2} Y_5 Y_6 Y_4 
  \underline{\Phi_4} \Phi_6^{-1} \Phi_2^{-1} \Phi_4^{-1}
  \tag*{by \eqref{theta3} and \eqref{eq:center}}  
  \\
  &=
  q \, \Phi_4 \underline{\Phi_2 \Phi_6} 
  Y_2 Y_5 Y_6 \underline{Y_4 (1+q Y_4)^{-1} 
  \Phi_6^{-1} \Phi_2^{-1}} \Phi_4^{-1}
  \tag*{by \eqref{R_difference}}  
  \\
  &=
  q \, \Phi_4 Y_2 Y_5 Y_6 Y_4 (1+q^{-1}Y_2)^{-1} (1+q^{-1}Y_6)^{-1}
  \left(1+q Y_4 (1+q^{-1}Y_2)^{-1} (1+q^{-1}Y_6)^{-1} \right)^{-1}
  \Phi_4^{-1}
  \tag*{by \eqref{R_difference}}  
  \\
  &=
  q^{-1} \Phi_4 Y_5 Y_4 
  (1+q Y_2^{-1})^{-1} (1+q Y_6^{-1})^{-1}
  \left(1+q Y_4 (1+q^{-1}Y_2)^{-1} (1+q^{-1}Y_6)^{-1} \right)^{-1}
  \Phi_4^{-1}
  \\
  &=
  Y_5 (1+q Y_4^{-1})^{-1}(1+q {Y_2^\prime}^{-1})^{-1}
  (1+q {Y_6^\prime}^{-1})^{-1} (1+q {Y_4^{\prime \prime}}^{-1})^{-1}
  \tag*{by \eqref{R_difference}}.
\end{align*}
Here we have used $\E^{2 \pi b c} = q Y_5 Y_6$ at the second equality,
and $Y_2^\prime$, $Y_6^\prime$ and $Y_4^{\prime \prime}$
are given by~\eqref{eq:Y246}.
A case of $i=3$ is as follows.
\begin{align*}
  \overset{1}{\mathcal{R}\vphantom{{}^{-1}}}  \, Y_3 \,
  \overset{1\phantom{-1}}{\mathcal{R}^{-1}}
  &=
  \Phi_4 \Phi_2 \Phi_6 \Phi_4^{-1} 
  \underline{\theta_{c+4} Y_3 \theta_{c+4}^{-1}} 
  \Phi_4 \Phi_6^{-1} \Phi_2^{-1} \Phi_4^{-1}
  \\
  &=
  q \, \Phi_4 \Phi_2 \Phi_6 \underline{\Phi_4^{-1} 
  Y_2^{-1} Y_4^{-1} \Phi_4} \Phi_6^{-1} \Phi_2^{-1} \Phi_4^{-1}
  \tag*{by \eqref{theta2} and \eqref{eq:center}}  
  \\
  &=
  \Phi_4 \underline{\Phi_2 \Phi_6}  
  Y_2^{-1} \underline{(1+q Y_4^{-1})
  \Phi_6^{-1} \Phi_2^{-1}} \Phi_4^{-1}
  \tag*{by \eqref{R_difference}}  
  \\
  &=    
  \Phi_4 
  Y_2^{-1} \left(1+q Y_4^{-1} (1+qY_2) (1+qY_6)\right)
  \Phi_4^{-1}
  \tag*{by \eqref{R_difference}}  
  \\
  &=
  {Y_2^{\prime}}^{-1} (1+q Y_4^{\prime\prime})
  \tag*{by \eqref{R_difference}},
\end{align*}
\endgroup
where we have used $\E^{-2 \pi b c} = q Y_2^{-1} Y_3^{-1}$ 
at the second equality.

\section{Proof of Braid Relation~\eqref{eq:braid-1}}
\label{sec:proof_YBE}

We shall check~\eqref{eq:braid-1} for $i=1$. 
We employ the notations in Appendix~\ref{sec:proof_adR}.
The proof is straightforward but tedious by use
of~\eqref{pentagon},~\eqref{theta2}, and~\eqref{theta3}.
We compute  as follows;
\begingroup
\allowdisplaybreaks
\begin{align*}
  \overset{2}{\mathcal{R}} \, \overset{1}{\mathcal{R}} 
  & =
  \Phi_7 \Phi_5 \Phi_9 \Phi_7^{-1} 
  \underline{\theta_{c+7}}
  \Phi_4 \Phi_2
  \underline{\Phi_6}
  \Phi_4^{-1} \theta_{c+4}
  \\
  & =
  \Phi_7 \Phi_5 \Phi_9 
  \underline{\Phi_7^{-1}}
  \Phi_4 \Phi_2
  \underline{\Phi_{6-7-c}}
  \Phi_4^{-1} \theta_{4+c} \theta_{7+c}
  \tag*{by \eqref{theta2}}
  \\
  & =
  \Phi_7 \Phi_5 
  \underline{\Phi_9}
  \Phi_4 \Phi_2 \Phi_{-5} 
  \underline{\Phi_{-5-7}}
  \Phi_4^{-1} \theta_{4+c} \Phi_7^{-1} \theta_{7+c}
  \tag*{by \eqref{pentagon}}
  \\
  & =
  \Phi_7 
  \underline{\Phi_5 \Phi_4} 
  \Phi_2 \Phi_{-5} \Phi_{-5-7} \Phi_{9-7-5} \Phi_4^{-1}
  \theta_{4+c} \Phi_9 \Phi_7^{-1} \theta_{7+c}
  \tag*{by \eqref{pentagon}}
  \\
  & =
  \Phi_7 \Phi_4 \Phi_{4+5} \Phi_2
  \underline{\Phi_{-5} \Phi_5} \, 
  \underline{\Phi_{-5-7} \Phi_{9-7-5}}
  \Phi_4^{-1} \theta_{4+c} \Phi_9 \Phi_7^{-1} \theta_{7+c}
  \tag*{by \eqref{pentagon}}
  \\
  & =
  \underline{\Phi_7} 
  \Phi_4
  \underline{\Phi_{4+5} }
  \Phi_2 \Phi_{-7} \Phi_{9-7} \Phi_{4+5}^{-1} \theta_{4+5+c}
  \Phi_{-5} \Phi_5 \Phi_9 \Phi_7^{-1} \theta_{7+c}
  \tag*{by \eqref{theta3}}
  \\
  & =
  \Phi_4 \Phi_{4+5} \Phi_{4+5+7} \Phi_2 
  \underline{\theta_7
    \Phi_{9-7} \Phi_{4+5}^{-1} \theta_{4+5+c} \Phi_{-5}
  }
  \Phi_5 \Phi_9 \Phi_7^{-1} \theta_{7+c}
  \tag*{by \eqref{pentagon}}
  \\
  & =
  \Phi_4 \Phi_{4+5} \Phi_{4+5+7} \Phi_2 \Phi_9 \Phi_{4+5+7}^{-1}
  \theta_{4+5+7+c} \Phi_{-5-7} \Phi_{-7}
  \overset{2}{\mathcal{R}}
  \tag*{by \eqref{theta3}}
\end{align*}
Thus we get
\begin{align*}
  \overset{1}{\mathcal{R}} \,
  \overset{2}{\mathcal{R}} \,
  \overset{1}{\mathcal{R}}
  & =
  \Phi_4 \Phi_2 \Phi_6 \Phi_4^{-1}
  \underline{\theta_{4+c}}
  \Phi_4
  \underline{\Phi_{4+5} \Phi_{4+5+7} \Phi_2 }
  \Phi_9 
  \underline{\Phi_{4+5+7}^{-1} \theta_{4+5+7+c} \Phi_{-5-7}}
  \Phi_{-7} 
  \overset{2}{\mathcal{R}}
  \\
  & =
  \Phi_4 \Phi_2 
  \underline{\Phi_6  \Phi_{5-c}} \,
  \underline{\Phi_{5+7-c}}
  \Phi_{2+4+c} \Phi_9 \Phi_{5+7-c}^{-1} \theta_{5+7} \Phi_{-5+4+c-7}
  \Phi_{-7} \theta_{4+c} 
  \overset{2}{\mathcal{R}}
  \tag*{by \eqref{theta2}}
  \\
  & =
  \Phi_4 \Phi_2 
  \underline{\theta_6}
  \Phi_{7-6}
  \Phi_{2+4+c} \Phi_9 \Phi_{5+7-c}^{-1} \theta_{5+7} \Phi_{-5+4+c-7}
  \Phi_{-7} \theta_{4+c} 
  \overset{2}{\mathcal{R}}
  \\
  & =
  \Phi_4
\Phi_2
  \Phi_7 \Phi_{2+4+5} \Phi_9 \Phi_7^{-1}
  \theta_{7+c} \Phi_{4-6-7} \Phi_{-6-7} \Phi_{4-6}
  \Phi_{-6} \Phi_6 \Phi_4^{-1} \theta_{4+c} 
  \overset{2}{\mathcal{R}}
  \tag*{by \eqref{theta3}}
  \\
  & =
  \Phi_4 \Phi_2 \Phi_7 \Phi_{2+4+5} \Phi_9 \Phi_7^{-1}
  \theta_{7+c} \Phi_{4-6-7} \Phi_{-6-7} \Phi_{4-6} \Phi_{-6}
  \Phi_2^{-1} \Phi_4^{-1}
  \overset{1}{\mathcal{R}}
  \overset{2}{\mathcal{R}}
  \tag*{by \eqref{pentagon}}
  \end{align*}
In the last expression
we have
\begin{align*}
  &
  \Phi_4 \Phi_2 \Phi_7 \Phi_{2+4+5} \Phi_9 \Phi_7^{-1}
  \theta_{7+c} \Phi_{4-6-7} \Phi_{-6-7} \Phi_{4-6} \Phi_{-6}
  \underline{\Phi_2^{-1} \Phi_4^{-1}}
  \\
  & =
  \Phi_4 \Phi_2 \Phi_7 \Phi_{2+4+5} \Phi_9 \Phi_7^{-1}
  \theta_{7+c} \Phi_{4-6-7} \Phi_{-6-7} 
\underline{\Phi_{4-6} \Phi_{-6} \Phi_4^{-1}}
\Phi_{2+4}^{-1} \Phi_2^{-1}
  \tag*{by \eqref{pentagon}}
  \\
  &=
  \Phi_4 \Phi_2 \Phi_7 \Phi_{2+4+5} \Phi_9 \Phi_7^{-1}
  \underline{\theta_{7+c}}
  \Phi_{4-6-7} \Phi_{-6-7} 
  \Phi_{4}^{-1} \Phi_{-6}
  \Phi_{2+4}^{-1} \Phi_2^{-1}
  \tag*{by \eqref{pentagon}}
  \\
  & =
  \Phi_4 \Phi_2 \Phi_7 \Phi_{2+4+5} \Phi_9
\Phi_7^{-1}
\underline{\Phi_{4+5} \Phi_5 \Phi_{4}^{-1}}
\Phi_{5+7}
  \Phi_{2+4}^{-1} \Phi_2^{-1}
  \theta_{7+c}
  \tag*{by \eqref{theta2}}
  \\
  & =
  \Phi_4 \Phi_2 \Phi_7 \Phi_{2+4+5} \Phi_9
  \Phi_7^{-1}   \Phi_{4}^{-1}
  \underline{
    \Phi_5
    \Phi_{5+7} \Phi_7
  }
  \Phi_{2+4}^{-1} \Phi_2^{-1}
  \Phi_7^{-1}
  \theta_{7+c}
  \\
  & =
  \Phi_4 \Phi_2 \Phi_7
  \underline{\Phi_{2+4+5} }
  \Phi_{4}^{-1}
  \underline{
    \Phi_5    \Phi_{2+4}^{-1}
  }
  \Phi_2^{-1}
  \Phi_9
  \Phi_7^{-1}
  \theta_{7+c}
  \\
  & =
  \Phi_4 \Phi_2 
  \underline{
    \Phi_4^{-1} \Phi_{2+4}^{-1} \Phi_2^{-1}
  }
  \Phi_7\Phi_5 
  \Phi_9 \Phi_7^{-1} \theta_{7+c}
  \tag*{by \eqref{pentagon}}
  \\
  & =
  \overset{2}{\mathcal{R}}
  \tag*{by \eqref{pentagon}}
\end{align*}
\endgroup
This completes the proof.

\section{\mathversion{bold}Proof of \eqref{Fourier_w}}
\label{sec:Fourier_w}
We recall a transformation of terminating $q$-hypergeometric function~\cite[(III.6)]{GaspRahm04}
\begin{equation*}
  \sum_{k=0}^n
  \frac{(q^{-n};q)_k \, (b;q)_k}{
    (c;q)_k \, (q;q)_k
  } \, z^k
  =
  \frac{
    (c/b;q)_n}{
    (c;q)_n}
  \, \left( \frac{b \, z}{q} \right)^n
  \sum_{k=0}^n
  \frac{
    (q^{-n};q)_k (q/z;q)_k (q^{1-n}/c;q)_k}{
    (b \, q^{1-n}/c; q)_k (q;q)_k
  }\, q^k ,
\end{equation*}
which reduces, in $c\to 0$, to
\begin{equation}
  \sum_{k=0}^n 
  \frac{(q^{-n};q)_k (b;q)_k}{(q;q)_k} \, z^k
  =
  \left( \frac{b \, z}{q} \right)^n
  \sum_{k=0}^n
  \frac{
    (q^{-n};q)_k (q/z;q)_k
  }{
    (q;q)_k
  } \, \left( \frac{q}{b} \right)^k .
\end{equation}
By setting $n=N-1$ and $q=\omega$ where $\omega^N=1$, we get
\begin{equation}
  \label{omega_identity}
  \sum_{k=0}^{N-1} (b;\omega)_k z^k
  =
  \left(
    \frac{b \, z}{\omega} 
  \right)^{N-1} 
  \sum_{k=0}^{N-1} (\omega/z; \omega)_k \,
  \left( \frac{\omega}{b} \right)^k .
\end{equation}

Using this identity, we compute as follows.
\begin{align*}
  \sum_{k=0}^{N-1} \frac{\omega^{-n k}}{
    w(x,y|k)
  }
  &=
  \frac{1}{
    w(x,y|0)}
  \sum_{k=0}^{N-1} ( \omega \, x; \omega)_k 
  \left(
    \frac{1}{y \, \omega^n }
  \right)^k
  \\
  & =
  \frac{1}{
    w(x,y|0)}
  \left( \frac{x}{y} \right)^{N-1} \omega^n
  \sum_{k=0}^{N-1} (y \, \omega^{n+1};\omega)_k \, x^{-k}
  \tag*{by \eqref{omega_identity}}
  \\
  & =
  \frac{w(y \, \omega^n, x|0)}{
    w(x,y|0)
  } \,
  \left( \frac{x}{y} \right)^{N-1} \omega^n \sum_{k=0}^{N-1}
  \frac{1}{w(y \, \omega^n, x|k )} .
\end{align*}
Here we note that $x^N+y^N=1$.
As a result, we get from~\eqref{w_modulo} and~\eqref{w_omega_n}
\begin{equation}
  \label{w_inverse}
  \sum_{k=0}^{N-1} \frac{\omega^{-n k}}{w(x,y|k)}
  =
  \omega^n \, w(y  \, \omega^n, x | 0) \, \left( \frac{x}{y}
  \right)^{\frac{N-1}{2}} \,
  \lambda(x,y),
\end{equation}
where
\begin{equation}
  \lambda(x,y)
  =
  \frac{
    (x/y)^{(N-1)/2}
  }{
    w(x,y|0)
  }
  \sum_{k=0}^{N-1} \frac{1}{w(y,x|k)} .
\end{equation}
See that 
\begin{gather}
  \lambda(x,\omega \, y) = \lambda(x, y),
  \\
  \lambda(x,y) \stackrel{x\to 1}{\to}
  N \, \prod_{j=1}^{N-1} (1 - \omega^{-j})^{-j/N} .
\end{gather}
The Fourier transform
of~\eqref{w_inverse} reduces to
\begin{equation}
  \sum_{k=0}^{N-1} w(y,x | k ) \, \omega^{k n}
  =
  N \,
  \frac{(y/x)^{(N-1)/2}}{\lambda(x,y)} \,
  \frac{1}{w(x,y|n-1)} .
\end{equation}

\end{document}